\documentclass[12pt,reqno,a4wide]{amsart}

\usepackage{tikz} 
\usepackage{calrsfs}
\usepackage[charter]{mathdesign}

\allowdisplaybreaks

\title{ Deutsch paths and their enumeration}

\author{Helmut Prodinger}
\address{Helmut Prodinger\\
	Department of Mathematical Sciences\\
	Stellenbosch University\\
	7602 Stellenbosch\\
	South Africa}
\email{hproding@sun.ac.za}

\begin{document}
\begin{abstract}
A variation of Dyck paths allows for down-steps of arbitrary length, not just one. Credits for this invention are given to Emeric Deutsch.
Surprisingly, the enumeration of them is somewhat akin to the analysis of Motzkin-paths; the last section contains a bijection.
	\end{abstract}

\maketitle
	
	\section{Introduction}
	
	Cameron~\cite{Cameron00} considers non-negative lattice paths with an up-step of one unit and down-steps of $2,4,6,\dots$ units.
	She gives (loc. cit.) credits to Emeric Deutsch for inventing these paths. Before we engage into further analysis of 
	this situation (in a future publication), we will consider here the more modest model of down-steps of $1,2,3,4,\dots$ units.
	Lattice paths with various step sizes have been thoroughly studied by Banderier and Flajolet~\cite{Banderier}, but  never for an infinite number of 	possible steps.
	
	We will call these paths \emph{Deutsch} paths,\footnote{The fact that \emph{Dyck paths} and \emph{Deutsch paths} sound so similar is also a bonus.} and require that they are non-negative and finish on a prescribed level $i$. If $i=0$, which is the analogue of Dyck paths, we talk about a \emph{closed Deutsch path}. Since the up-steps and down-steps are not symmetrical,
	we can also investigate the model where the paths develop from right to left. In this case we talk about \emph{reversed Deutsch paths}.

	Somewhat surprisingly, the computations that we will perform, lead us the \emph{world of Motzkin} paths.
	Recall that the generating function of Motzkin paths (up-step, level-step, down-step) leads to $M=1+zM+z^2M^2$,
	or
	\begin{equation*}
M(z)=\frac{1-z-\sqrt{1-2z-3z^2}}{2z^2}=1+z+2\,{z}^{2}+4\,{z}^{3}+9\,{z}^{4}+21\,{z}^{5}+51\,{z}^{6}+\cdots
	\end{equation*}
Here, as first shown in \cite{Prodinger-three}, the substitution $z=\dfrac{v}{1+v+v^2}$ is beneficial:
\begin{equation*}
M=1+v+v^2,\quad\text{and}\quad v=\frac{1-z-\sqrt{1-2z-3z^2}}{2z}.
\end{equation*}
	Introducing \emph{trinomial coefficients} via
	\begin{equation*}
\binom{n,3}{k}:=[v^k](1+v+v^2)^n
	\end{equation*}
	(notation is from Comtet~\cite{Comtet-book}), we can compute the Motzkin numbers
\begin{align*}
	[z^n]M(z)&=\frac1{2\pi i}\oint \frac{dz}{z^{n+1}}M(z)\\*
	&=\frac1{2\pi i}\oint \frac{dv}{v^{n+1}}(1-v^2)(1+v+v^2)^{n}\\
	&=\binom{n,3}{n}-\binom{n,3}{n-2}.
\end{align*}	
	
	For the analysis of Deutsch paths, we will see that the same substitution works extremely well, and explicit answers can be given in terms of trinomial coefficients. The last section contains a bijection between open Deutsch paths (level at the end arbitrary) and Motzkin paths, of the same length.

	\section{Enumeration of Deutsch paths}
	We introduce an upper boundary $h$, so that the paths live  in a strip. This has the advantage that generating functions can be described by a system of linear equations, and solving it can be done by Cramer's rule. This involves the compution of some determinants. 
	
	We present the example of $h=4$, and the generating functions $\varphi_i(z)$ describe bounded Deutsch paths ending at level $i$.
	
	\begin{equation*}
		\left(\begin{matrix}
			1&-z&-z&-z&-z\\
			-z&	1&-z&-z&-z\\
			0&	-z&	1&-z&-z\\
			0& 0&	-z&	1&-z\\
			0& 0& 0&	-z&	1\\
		\end{matrix}\right)
		\left(\begin{matrix}
			\varphi_0\\
			\varphi_1\\
			\varphi_2\\
			\varphi_3\\
			\varphi_4\\
		\end{matrix}\right)=
		\left(\begin{matrix}
			1\\
			0\\
			0\\
			0\\
			0\\
		\end{matrix}\right)
	\end{equation*}	
The 	$n\times n $ determinant
	\begin{equation*}
D_n=\frac{(1+v)^{n-1}}{(1+v+v^2)^n}\frac{1-v^{n+2}}{1-v}
	\end{equation*}
	of this system is of importance. It will eventually turn out that
\begin{equation*}
	D_{h+1}=\frac{(1+v)^{h}}{(1+v+v^2)^{h+1}}\frac{1-v^{h+3}}{1-v},
\end{equation*}
which can be obtained from the recursion
\begin{equation*}
	(1+v+v^2)^2D_{n+2}-(1+v+v^2)(1+v)^2D_{n+1}+v(1+v)^2D_{n}=0.
\end{equation*}	
Applying Cramer's rule leads then to
	\begin{align*}
\varphi_i&=\frac{z^iD_{h-i}}{D_{h+1}}=\frac{v^i}{(1+v+v^2)^{i}}
\frac{(1+v+v^2)^{h+1}}{(1+v)^{h}}\frac{1-v}{1-v^{h+3}}\frac{(1+v)^{h-i-1}}{(1+v+v^2)^{h-i}}\frac{1-v^{h-i+2}}{1-v}\\
&=\frac{v^i(1+v+v^2)}{(1+v)^{i+1} }\frac{1-v^{h-i+2}}{1-v^{h+3}}.
	\end{align*}

	Of special interest is the generating function of closed bounded Deutsch paths:
	\begin{align*}
		\varphi_0
		&=\frac{1+v+v^2}{1+v }\frac{1-v^{h+2}}{1-v^{h+3}}.
	\end{align*}
	Taking the limit $h\to\infty$ leads to the enumeration of closed Deutsch paths, as originally desired:
	\begin{align*}
		\varphi_0
		&=\frac{1+v+v^2}{1+v }.
	\end{align*}
We can read off the coefficients explicitly:
\begin{align*}
	[z^n]\varphi_0&=\frac1{2\pi i}\oint \frac{dz}{z^{n+1}}\frac{1+v+v^2}{1+v }\\
	&=\frac1{2\pi i}\oint \frac{dv}{v^{n+1}}(1-v)(1+v+v^2)^{n}\\
	&=\binom{n,3}{n}-\binom{n,3}{n-1}.
\end{align*}	
	
The enumeration of closed Deutsch paths  of height $\ge h$
	\begin{align*}
\frac{1+v+v^2}{1+v }&-\frac{1+v+v^2}{1+v }\frac{1-v^{h+1}}{1-v^{h+2}}
=\frac{1+v+v^2}{1+v }\frac{v^{h+1}(1-v)}{1-v^{h+2}}\\
	\end{align*}
is of interest, and we will discuss asymptotics of the average height in a later section.

Here is the generating function of bounded Deutsch paths with (arbitrary) open end:
\begin{align*}
	\sum_{i=0}^h\varphi_i
	&=(1+v+v^2)\frac{1-v^{h+1}}{1-v^{h+3}}.
\end{align*}
In the limit $h\to\infty$, this leads to
\begin{align*}
1+v+v^2,
\end{align*}
which is the generating function of Motzkin paths, according to length. A combinatorial explanation is given in the last section.

\section{Enumeration of reversed Deutsch paths}

This enumeration is quite similar, only the matrix is transposed:
			\begin{equation*}
			\left(\begin{matrix}
				1&-z&0&0&0\\
				-z&	1&-z&0&0\\
				-z&	-z&	1&-z&0\\
				-z& -z&	-z&	1&-z\\
				-z& -z& -z&	-z&	1\\
			\end{matrix}\right)
			\left(\begin{matrix}
				\psi_0\\
				\psi_1\\
				\psi_2\\
				\psi_3\\
				\psi_4\\
			\end{matrix}\right)=
			\left(\begin{matrix}
				1\\
				0\\
				0\\
				0\\
				0\\
			\end{matrix}\right)
		\end{equation*}	
		The determinant is the same as before, but the application of Cramer's rule is slightly more unpleasant, as we have to distinguish cases. We only present the final results; intermediate calculations are not really difficult, but require some concentration.
	\begin{equation*}
\psi_0=\frac{1+v+v^2}{1+v}\frac{1-v^{h+2}}{1-v^{h+3}}
	\end{equation*}
	and for $i\ge1$:
	\begin{equation*}
		\psi_i=v(1+v+v^2)(1+v)^{i-2}\frac{1-v^{h+1-i}}{1-v^{h+3}}
	\end{equation*}
		\begin{equation*}
\sum_{i=0}^h\psi_i=(1+v+v^2)(1+v)^h\frac{1-v}{1-v^{h+3}},
	\end{equation*}
and the limit of this for $h\to\infty$:
	\begin{equation*}
		(1+v+v^2)(1-v).
	\end{equation*}
The enumeration in terms of trinomial coefficients:
	\begin{equation*}
\binom{n,3}{n}-\binom{n,3}{n-1}-\binom{n,3}{n-2}+\binom{n,3}{n-3}.
	\end{equation*}

	\begin{figure}	
		\begin{center}
			\begin{tikzpicture}[scale=0.6]
			
			\draw[step=1.cm,black,dotted] (-0.0,-0.0) grid (3.0,1.0);

			\draw[ultra thick] (0,0) to (1,1) to (2,0) to (3,0);
			
			\end{tikzpicture}
			\begin{tikzpicture}[scale=0.6]
			
			\draw[step=1.cm,black,dotted] (-0.0,-0.0) grid (3.0,1.0);

			\draw[ultra thick] (0,0) to (1,1) to (2,1) to (3,0);
			
			\end{tikzpicture}
			\begin{tikzpicture}[scale=0.6]
			
			\draw[step=1.cm,black,dotted] (-0.0,-0.0) grid (3.0,1.0);

			\draw[ultra thick] (0,0) to (1,0) to (2,1) to (3,0);
			
			\end{tikzpicture}
			\begin{tikzpicture}[scale=0.6]
			
			\draw[step=1.cm,black,dotted] (-0.0,-0.0) grid (3.0,1.0);

			\draw[ultra thick] (0,0) to (1,0) to (2,0) to (3,0);
			
			\end{tikzpicture}

			\end	{center}
			\caption{4 Motzkin paths of length 3}
		\end{figure}
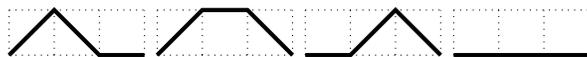
		
		\begin{figure}	
			\begin{center}
				\begin{tikzpicture}[scale=0.6]
				
				\draw[step=1.cm,black,dotted] (-0.0,-0.0) grid (3.0,3.0);

				\draw[ultra thick] (0,0) to (1,1) to (2,0) to (3,1);
				
				\end{tikzpicture}
				\begin{tikzpicture}[scale=0.6]
				
				\draw[step=1.cm,black,dotted] (-0.0,-0.0) grid (3.0,3.0);

				\draw[ultra thick] (0,0) to (1,1) to (2,2) to (3,1);
				
				\end{tikzpicture}
				\begin{tikzpicture}[scale=0.6]
				
				\draw[step=1.cm,black,dotted] (-0.0,-0.0) grid (3.0,3.0);
				
				\draw[ultra thick] (0,0) to (1,1) to (2,2) to (3,0);
				
				\end{tikzpicture}
				\begin{tikzpicture}[scale=0.6]
				
				\draw[step=1.cm,black,dotted] (-0.0,-0.0) grid (3.0,3.0);

				\draw[ultra thick] (0,0) to (1,1) to (2,2) to (3,3);
				
				\end{tikzpicture}

				\end	{center}
				\caption{4 open Deutsch paths   of length 3}
			\end{figure}
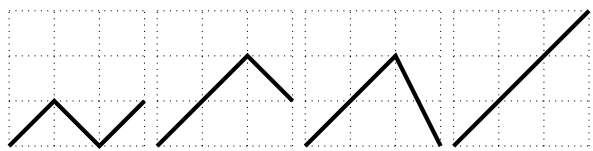
			
			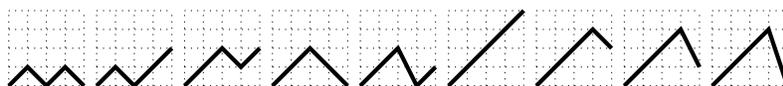
\begin{figure}	
				\begin{center}
					\begin{tikzpicture}[scale=0.25]
					\draw[step=1.cm,black,dotted] (-0.0,-0.0) grid (4.0,4.0);
					\draw[ultra thick] (0,0) to (1,1) to (2,0) to (3,1) to (4,0);
					\end{tikzpicture}
					\begin{tikzpicture}[scale=0.25]
					\draw[step=1.cm,black,dotted] (-0.0,-0.0) grid (4.0,4.0);
					\draw[ultra thick] (0,0) to (1,1) to (2,0) to (3,1)to (4,2);
					\end{tikzpicture}
					\begin{tikzpicture}[scale=0.25]
					\draw[step=1.cm,black,dotted] (-0.0,-0.0) grid (4.0,4.0);
					\draw[ultra thick] (0,0) to (1,1) to (2,2) to (3,1)to(4,2);
					\end{tikzpicture}
					\begin{tikzpicture}[scale=0.25]
					\draw[step=1.cm,black,dotted] (-0.0,-0.0) grid (4.0,4.0);
					\draw[ultra thick] (0,0) to (1,1) to (2,2) to (3,1)to (4,0);
					\end{tikzpicture}
					\begin{tikzpicture}[scale=0.25]
					\draw[step=1.cm,black,dotted] (-0.0,-0.0) grid (4.0,4.0);
					\draw[ultra thick] (0,0) to (1,1) to (2,2) to (3,0)to(4,1);
					\end{tikzpicture}
					\begin{tikzpicture}[scale=0.25]
					\draw[step=1.cm,black,dotted] (-0.0,-0.0) grid (4.0,4.0);
					\draw[ultra thick] (0,0) to (1,1) to (2,2) to (3,3)to(4,4);
					\end{tikzpicture}
					\begin{tikzpicture}[scale=0.25]
					\draw[step=1.cm,black,dotted] (-0.0,-0.0) grid (4.0,4.0);
					\draw[ultra thick] (0,0) to (1,1) to (2,2) to (3,3)to(4,2);
					\end{tikzpicture}
					\begin{tikzpicture}[scale=0.25]
					\draw[step=1.cm,black,dotted] (-0.0,-0.0) grid (4.0,4.0);
					\draw[ultra thick] (0,0) to (1,1) to (2,2) to (3,3)to(4,1);
					\end{tikzpicture}
					\begin{tikzpicture}[scale=0.25]
					\draw[step=1.cm,black,dotted] (-0.0,-0.0) grid (4.0,4.0);
					\draw[ultra thick] (0,0) to (1,1) to (2,2) to (3,3)to(4,0);
					\end{tikzpicture}
					\end	{center}
					\caption{9 open Deutsch paths   of length 4}
				\end{figure}
				
				\begin{figure}	
					\begin{center}
						\begin{tikzpicture}[scale=0.25]
						\draw[step=1.cm,black,dotted] (-0.0,-0.0) grid (4.0,2.0);
						\draw[ultra thick] (0,0) to (1,1) to (2,0) to (3,1) to (4,0);
						\end{tikzpicture}
						\begin{tikzpicture}[scale=0.25]
						\draw[step=1.cm,black,dotted] (-0.0,-0.0) grid (4.0,2.0);
						\draw[ultra thick] (0,0) to (1,1) to (2,2) to (3,1)to (4,0);
						\end{tikzpicture}
						\begin{tikzpicture}[scale=0.25]
						\draw[step=1.cm,black,dotted] (-0.0,-0.0) grid (4.0,2.0);
						\draw[ultra thick] (0,0) to (1,0) to (2,0) to (3,1)to (4,0);
						\end{tikzpicture}
						\begin{tikzpicture}[scale=0.25]
						\draw[step=1.cm,black,dotted] (-0.0,-0.0) grid (4.0,2.0);
						\draw[ultra thick] (0,0) to (1,0) to (2,1) to (3,1)to (4,0);
						\end{tikzpicture}
						\begin{tikzpicture}[scale=0.25]
						\draw[step=1.cm,black,dotted] (-0.0,-0.0) grid (4.0,2.0);
						\draw[ultra thick] (0,0) to (1,0) to (2,1) to (3,0)to (4,0);
						\end{tikzpicture}
						\begin{tikzpicture}[scale=0.25]
						\draw[step=1.cm,black,dotted] (-0.0,-0.0) grid (4.0,2.0);
						\draw[ultra thick] (0,0) to (1,0) to (2,0) to (3,0)to (4,0);
						\end{tikzpicture}
						\begin{tikzpicture}[scale=0.25]
						\draw[step=1.cm,black,dotted] (-0.0,-0.0) grid (4.0,2.0);
						\draw[ultra thick] (0,0) to (1,1) to (2,1) to (3,1)to (4,0);
						\end{tikzpicture}
						\begin{tikzpicture}[scale=0.25]
						\draw[step=1.cm,black,dotted] (-0.0,-0.0) grid (4.0,2.0);
						\draw[ultra thick] (0,0) to (1,1) to (2,1) to (3,0)to (4,0);
						\end{tikzpicture}
						\begin{tikzpicture}[scale=0.25]
						\draw[step=1.cm,black,dotted] (-0.0,-0.0) grid (4.0,2.0);
						\draw[ultra thick] (0,0) to (1,1) to (2,0) to (3,0)to (4,0);
						\end{tikzpicture}

						\end	{center}
						\caption{9 Motzkin paths   of length 4}
					\end{figure}
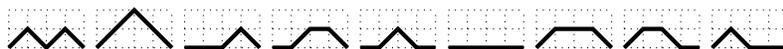
				
	\section{The average height of closed and open Deutsch paths}				
	
	It was already worked out what the generating function of closed Deutsch paths is, and for the average height, one needs to compute	
	\begin{align*}
		\sum_{h\ge1}&\frac{1+v+v^2}{1+v }\frac{v^{h+1}(1-v)}{1-v^{h+2}}\\
		&=-\frac{1+v+v^2}{1+v }
		-\frac{v(1+v+v^2)}{(1+v)^2 }+	\frac{(1+v+v^2)(1-v)}{v(1+v) }\sum_{h\ge1}\frac{v^{h}}{1-v^{h}}\\
	\end{align*}
	
The standard procedure is to find the local behaviour of this near the singularity $v=1\leftrightarrow z=\frac13$. A fairly detailed description of this can be found in \cite{Heuberger}. Somewhat similar computations will also appear in \cite{Baril}. We get
	\begin{align*}
	E(z)&=	\sum_{h\ge1}\frac{1+v+v^2}{1+v }\frac{v^{h+1}(1-v)}{1-v^{h+2}}\\
		&=-\frac94 +	\frac{3(1-v)}{2}\Big[-\frac{\log(1-v)}{1-v}+\frac{\gamma}{1-v}+\cdots\Big]\\
		&=-\frac32\log(1-v)+\textsf{const.}+\cdots\\
	\end{align*}
Since $1-v\sim\sqrt3\sqrt{1-3z}$, this leads further to the local expansion
\begin{equation*}
E(z)\sim-\frac34\log(1-3z)+\cdots
\end{equation*}	
Traditional singularity analysis~\cite{FlOd} 		leads to
\begin{equation*}
[z^n]E(z)\sim
\frac34\frac{3^n}{n}.
\end{equation*}
This has to be divided by the total number of closed Deutsch paths, with generating function
\begin{equation*}
\frac{1+v+v^2}{1+v}\sim 3-\frac{9}{4\sqrt3}\sqrt{1-3z},
\end{equation*}					
so that
\begin{equation*}
[z^n]\frac{1+v+v^2}{1+v}\sim \frac{9}{8\sqrt{3\pi}}3^nn^{-3/2}.
\end{equation*}
The quotients is the average height of closed Deutsch paths (leading term only):
\begin{equation*}
\frac{\frac34\frac{3^n}{n}}{\frac{9}{8\sqrt{3\pi}}3^nn^{-3/2}}=2\sqrt{\frac{\pi n}{3}}.
\end{equation*}
In the classical case of Motzkin paths \cite{Prodinger-three}, this height is $\sqrt{\frac{\pi n}{3}}$, so Deutsch paths are about twice as high.

We would also like to do this type of analysis for open Deutsch paths, as they are exactly enumerated by Motzkin numbers.
Recall the generating function
\begin{align*}
	\sum_{i=0}^h\varphi_i
	&=(1+v+v^2)\frac{1-v^{h+1}}{1-v^{h+3}}\\
\end{align*}
of open Deutsch paths, with height $\le h$. Taking differences, we get the generating function of open Deutsch paths, with height $\ge h$:
\begin{equation*}
(1+v+v^2)(1-v^2)\frac{v^h}{1-v^{h+2}}.
\end{equation*}
This has to be summed:
\begin{equation*}
(1+v+v^2)(1-v^2)\sum_{h\ge1}\frac{v^h}{1-v^{h+2}}\sim -6\log(1-v)\sim-3\log(1-3z),
\end{equation*}
with coefficients asymptotic to $\frac{3^{n+1}}{n}$.

This has to be divided by 
\begin{equation*}
[z^n]M(z)\sim\frac{9}{2\sqrt{3\pi }}3^nn^{-3/2};
\end{equation*}
which leads again to the average height 
\begin{equation*}
2\sqrt{\frac{\pi n}{3}}.
\end{equation*}
So open versus closed does not influence the (leading term of the) average height.

	\section{The LU-decomposition}
	
	We compute the
	LU-decomposition of the Deutsch matrix, where indices run from 1 to $n$. This isn't strictly necessary, but always interesting, 
	and leads to the determinant as a by-product.
	\begin{equation*}
L_{i,i}=1,
	\end{equation*}
	\begin{equation*}
		L_{i,i-1}=-\frac{v}{1+v}\frac{1-v^i}{1-v^{i+1}},
	\end{equation*}
	other values 0.
	
	\begin{equation*}
U_{i,i}=\frac{1+v}{1+v+v^2}\frac{1-v^{i+2}}{1-v^{i+1}},
	\end{equation*}
	\begin{equation*}
		U_{i,j}=-v\frac{1+v}{1+v+v^2}\frac{1-v^{i}}{1-v^{i+1}},
	\end{equation*}
	
	For the conformation, we make this computation:
	\begin{align*}
\sum_k L_{i,k}U_{k,j}&=U_{i,j}-\frac{v}{1+v}\frac{1-v^i}{1-v^{i+1}}U_{i-1,j}\\
&=-v\frac{1+v}{1+v+v^2}\frac{1-v^{i}}{1-v^{i+1}}+\frac{v}{1+v}\frac{1-v^i}{1-v^{i+1}}
v\frac{1+v}{1+v+v^2}\frac{1-v^{i-1}}{1-v^{i}}\\
&=-\frac{v}{1+v+v^2}=-z.
	\end{align*}
	Now let us consider the special cases, where $U_{i,i}$ appears:
	\begin{align*}
		\sum_k L_{i,k}U_{k,i}&=U_{i,i}-\frac{1+v}{1+v+v^2}\frac{1-v^{i+2}}{1-v^{i+1}}\\
		&=\frac{1+v}{1+v+v^2}\frac{1-v^{i+2}}{1-v^{i+1}}+
		\frac{v}{1+v}\frac{1-v^i}{1-v^{i+1}}
		v\frac{1+v}{1+v+v^2}\frac{1-v^{i-1}}{1-v^{i}}=1;
	\end{align*}
\begin{align*}
	\sum_k L_{i,k}U_{k,i-1}&=U_{i,i-1}-\frac{v}{1+v}\frac{1-v^i}{1-v^{i+1}}U_{i-1,i-1}\\
	&=-v\frac{1+v}{1+v+v^2}\frac{1-v^{i}}{1-v^{i+1}}
	-\frac{v}{1+v}\frac{1-v^i}{1-v^{i+1}}
	\frac{1+v}{1+v+v^2}\frac{1-v^{i+1}}{1-v^{i}}\\
	&=-\frac{v}{1+v+v^2}=-z.
\end{align*}
	In particular, we get for the determinant
	\begin{equation*}
U_{1,1}\dots U_{n,n}=\frac{(1+v)^n}{(1+v+v^2)^n}\frac{1-v^{n+1}}{1-v^{2}},
	\end{equation*}
	which checks with our previous observation.
	
	Very briefly, we only cite the results of the transposed matrix (related to reversed Deutsch paths):
		\begin{equation*}
L_{i,i}=1
	\end{equation*}
	and for $j<i$
	\begin{equation*}
		L_{i,j}=-v\frac{1-v^j}{1-v^{j+2}}.
	\end{equation*}
	Furthermore
	\begin{equation*}
		U_{i,i}=\frac{1+v}{1+v+v^2}\frac{1-v^{i+2}}{1-v^{i+1}},
	\end{equation*}
	\begin{equation*}
U_{i,i+1}=-\frac{v}{1+v+v^2},
	\end{equation*}
	and zero for all other values of the $U$-matrix.

	\section{The total area}

The area of a closed Deutsch path $a_0a_1\dots a_n$ with $a_0=a_n=0$ is defined to be $a_0+\cdots+a_n$.
With the generating functions that we worked out, the generating function of the total area, summed over all Deutsch
paths of length $n$, can be computed:
	
	\begin{align*}
A(z)&=\sum_{i\ge1}i\cdot\varphi_i\cdot\psi_i=\sum_{i\ge1}i\frac{v^i(1+v+v^2)}{(1+v)^{i+1}}v(1+v+v^2)(1+v)^{i-2}=
\frac{v^2(1+v+v^2)^2}{(1+v)^3(1-v)^2}\\
&={z}^{2}+3\,{z}^{3}+12\,{z}^{4}+39\,{z}^{5}+129\,{z}^{6}+411\,{z}^{7}+
1300\,{z}^{8}+4065\,{z}^{9}+12633\,{z}^{10}+\cdots
	\end{align*}
	We can also compute that
	\begin{equation*}
A(z)\sim \frac 98\frac1{(1-v)^2}\sim \frac38\frac{1}{1-3z},
	\end{equation*}
so that
\begin{equation*}
[z^n]A(z)\sim \frac38 3^n.
\end{equation*}
Dividing this by the total number of closed Deutsch paths
\begin{equation*}
\frac{9}{8\sqrt{3\pi }}3^nn^{-3/2}
\end{equation*}
leads to
\begin{equation*}
\sqrt{\frac{\pi}{3}}n^{3/2}.
\end{equation*}
One might even divide this by $n$, the length of the Deutsch paths, and can interpret $\sqrt{\frac{\pi n}{3}}$ as the average elevation of a random closed Deutsch path of length $n$. Notice that the average height is about twice this quantity.

	\section{A bijection}
	
	The generating function of Motzkin paths satisfies $M=1+zM+z^2M^2$, which is based on a first return decomposition. We will now derive the same recursion for open Deutsch paths, which leads to a recursively defined bijection.

We will define a map $\psi$	in a recursive way, mapping an open Deutsch path to a Motzkin path; it is easily seen to be reversible.
Of course the empty path is mapped to the empty path. Now, if the open Deutsch path never returns to the $x$-axis, it can be described by $w=U\widetilde{w}$, where $U$ describes an up-step, and $\widetilde{w}$ is itself an open Deutsch path. Then we map it to $F\psi(\widetilde{w})$, where $F$ is a flat step. Otherwise, if the open Deutsch path returns to the $x$-axis for the first time, it may be written as $w=U\widetilde{w}Dx$, where $D$ is a down-step of any size. Note that both, $\widetilde{w}$ and $x$ are themselves open Deutsch paths. Then we map this to $U\psi(\widetilde{w})D\psi(x)$, which is a Motzkin path (here, $D$ is just a down-step of one unit).

\newpage	
	
	\bibliographystyle{plain}

\end{document}